\documentclass[11pt,oneside,superscriptaddress,noshowpacs,centertags]{article}
\usepackage{afterpage}
\usepackage{fancyhdr}
\usepackage{etoolbox}
\usepackage{euscript}
\RequirePackage[square,numbers]{natbib}
\usepackage[colorlinks]{hyperref}
\usepackage{lipsum} %just to generate text for the example
\newcommand\shorttitle{Clairaut anti-invariant Riemannian maps to Sasakian manifolds}
\newcommand\authors{M. N. Zafar, A. Zaidi and G. Shanker}

\fancyheadoffset{\dimexpr\oddsidemargin+1in-\evensidemargin}

\usepackage{supertabular, setspace,hyperref}
\usepackage[a4paper, total={6.5in, 9.5in}]{geometry}
\usepackage{comment}
\usepackage{amsmath}
\usepackage{mathtools}
\usepackage{amsfonts}
\usepackage{amssymb}
\usepackage{amsthm}
\usepackage{enumitem}
\usepackage{caption}
\usepackage{subcaption}
\usepackage{authblk}
\newtheorem{theorem}{Theorem}[section]
\newtheorem{example}{Example}[section]
\newtheorem{definition}{Definition}[section]
\newtheorem{lemma}{Lemma}[section]

\usepackage{chngcntr}
\makeatletter
\patchcmd{\thebibliography}{\c@NAT@ctr\z@}{\c@NAT@ctr\z@
  \setlength{\labelwidth}{1.2em}%
  \setlength{\labelsep}{.5em}%
  \setlength{\leftmargin}{\dimexpr\labelwidth+\labelsep}%
}{}{}
\makeatother
\counterwithin*{equation}{section}
\counterwithin*{equation}{section}
\begin{document}
	\clearpage
	\date{}

	\title{\textbf {Clairaut anti-invariant Riemannian maps to Sasakian manifolds}}

	\author{Md Nadim Zafar}
	\author{Adeeba Zaidi}
	\author{Gauree Shanker\footnote{corresponding author}}
	\affil{\footnotesize Department of Mathematics and Statistics,
		Central University of Punjab, Bathinda, Punjab-151 401, India.\\
		Email: nadim.zafar13@gmail.com, adeebazaidi.az25@gmail.com and gauree.shanker@cup.edu.in*}
	\maketitle
	\begin{abstract}
		In this paper, we investigate the geometry of Clairaut anti-invariant Riemannnian maps whose base space are Sasakian manifolds. We obtain the necessary and sufficient conditions for a curve on a base manifold to be geodesic. We obtain conditions for an anti-invariant Riemannian map to be Clairaut. Further, we discuss the biharmonicity of such maps and construct some illustrative examples.
	\end{abstract}
	\vspace{0.2 cm}
 
\begin{small}
\end{small}
	
	
	\maketitle
	\section{Introduction}
	O'Neill\cite{O} and Gray\cite{G} were the first to describe Riemannian submersion between Riemannian manifolds. Later, in \cite{8}, a thorough discussion of the geometry of Riemannian submersions is provided. In 1992, Riemannian map between Riemannian manifolds was introduced by Fischer\cite{F} as a generalization of an isometric immersion and Riemannian submersion. 
	
	In addition, the geometry of the Riemannian map has been studied from different points of view and between various structures. \c{S}ahin introduced invariant and anti-invariant Riemannian maps to K{\"a}hler manifolds, semi-invariant Riemannian maps, slant Riemannian maps etc. (see \cite{S7}) taking the total space as an almost Hermitian manifold. Many authors studied Riemannian maps on contact manifolds such as cosymplectic, Sasakian, and Kenmotsu manifolds \cite{P1,P2,P3,P4,Azaidi}. In the book \cite{S7}, we find many recent developments in this field.
	
	In the theory of surfaces, Clairaut's theorem states that for any geodesic $\gamma$ on a surface $S$, the function $r\sin\theta$ is constant along $\gamma$, where $r$ is the distance from a point on the surface to the rotation axis and $\theta$ is the angle between $\gamma$ and the meridian through $\gamma$. In 1972, Bishop\cite{5} utilized this concept for Riemannian submersions and provided a necessary and sufficient condition for a Riemannian submersion to be a Clairaut-Riemannian submersion. Clairaut submersions were also investigated in \cite{AD,L,HMT}. In \cite{S11}, \c{S}ahin introduced the Clairaut Riemannian map using a geodesic curve on the total space and obtained the necessary and sufficient criteria for a Riemannian map to be a Clairaut Riemannian map. Later, Yadav and Meena introduced Clairaut invariant and anti-invatirant Riemannian maps from and to K\"ahler manifolds in \cite{AY,AY2,M2}. Recently, Li et al. \cite{YL} studied Clairaut semi-invariant Riemannian maps from cosymplectic manifolds. The investigation of diverse types of Riemannian maps has emerged as a prominent area of interest among researchers.
 	
	This paper presents the  geometry of Clairaut anti-invariant Riemannian map from a Riemannian manifold to a Sasakian manifold. The paper is structured as follows: In Section 2, we review the essential terminology and definitions used throughout the paper. In Section 3, we establish the Clairaut anti-invariant Riemannian map from a Riemannian manifold to a Sasakian manifold using a geodesic curve on the base space. We give the necessary and sufficient conditions for such a map to be a Clairaut map. Furthermore, we study the bi-harmonicity of these maps and construct some non-trivial examples.

	\section{Preliminaries}
	We begin this section by introducing the concepts of contact manifolds, Sasakian manifolds, and some key properties associated with Riemannian maps.
	A $(2n+1)$-dimensional smooth manifold $B$ is considered to have an almost-contact structure $(\varphi,\xi,\eta)$ if there exists a tensor field $\varphi$ of type $(1,1)$, a vector field $\xi$ (known as the characteristic or Reeb vector field), and a $1$-form $\eta$, which satisfy the following conditions \cite{B}:
	\begin{equation} \label{a}
		\varphi^2=-I + \eta \otimes \xi,~~~ \varphi \xi=0,~~\eta\circ\varphi=0,~~ \eta(\xi)=1.  
	\end{equation}
	Furthermore, an almost-contact manifold $B$ with a Riemannian metric $g_B$ is said to be compatible with the almost-contact structure $(\varphi, \xi, \eta)$ if it satisfies the following relations for any vector fields $B_1,B_2 \in \varGamma(TB)$ \cite{B}:
	\begin{align}
		& \ \ \ \ \ g_B(\varphi B_1,\varphi B_2)= g_B(B_1,B_2)-\eta (B_1)\eta (B_2),  \label{b} \\
		&g_B(\varphi B_1,B_2)= -g_B(B_1,\varphi B_2), \ \ \ \  \eta(B_1) = g_B(B_1,\xi) \label{c}
	\end{align}
	and the structure $(\varphi, \xi, \eta,g_B)$ is referred to as an almost contact metric structure. The almost-contact structure $(\varphi, \xi, \eta)$ is called normal if $\mathcal{N} + d\eta \otimes \xi=0$, where $\mathcal{N}$ is the Nijenhuis tensor of $\varphi$. Additionally, if $d\eta=\Phi$, where $\Phi(B_1,B_2)=g_B(\varphi B_1,B_2)$ is a tensor field of type $(0,2)$, then an almost contact metric structure is said to be a normal contact metric structure. A normal contact metric manifold $B$ is called a Sasakian manifold \cite{B} if it satisfies the following conditions:
	\begin{align}
		(\nabla_{B_{1}}\varphi)B_2 \ &= \ g_B(B_1,B_2)\xi-\eta(B_2)B_1,\label{d} \\  \nabla_{B_1}\xi \ &=-\varphi B_1, \label{e}
	\end{align}
	where $\nabla$ is the Levi-Civita connection of $g_B$ on $B$. Finally, the following example is a worth noting illustration of a Sasakian manifold on $\mathbb{R}^{2n+1}$.
	\begin{example}\cite{B}\label{S}
	Let us consider $\mathbb{R}^{2m+1}$ with Cartesian coordinates  $(x_i,y_i,z),\\i=1,...,m$ to be a Riemannian manifold with Riemannian metric $$g=\dfrac{1}{4}
		\begin{pmatrix}
			\delta _{ij}+y^iy^j & 0 & -y^i\\
			0 & \delta_{ij} & 0\\
			-y^i & 0  & 1 
		\end{pmatrix},$$  and  1-form $\eta=\frac{1}{2}(dz-\sum_{i=1}^{m}y^idx^i)$. The Reeb vector field is given by $2\frac{\partial}{\partial z}$, 
		and tensor field $\varphi$ is defined by
		$$\begin{pmatrix}
			0 & \delta _{ij} & 0 \\
			-\delta_{ij} & 0  & 0\\
			0 & y^i & 0  
		\end{pmatrix}.$$
		This structure gives a contact metric structure on $\mathbb{R}^{2m+1}$ and the vector fields $E_i=2\frac{\partial}{\partial y^i},\\E_{m+i}=2(\frac{\partial}{\partial x^i}+y^i\frac{\partial}{\partial z}), i=1,...,m$ and $\xi$ form a $\varphi$-basis for the contact metric structure. Also, it can be shown that $\mathbb{R}^{2m+1}$ is a Sasakian manifold with structure $(\varphi, \xi, \eta,g).$
  \end{example}			 	
	Further, consider a smooth map $F:(N^n,g_N)\rightarrow(B^b,g_B)$ between Riemannian manifolds $N$ and $B$. The differential map $F_*$ can be interpreted as a section of the bundle $Hom(TN,F^{-1}TB)\rightarrow N$, where $F^{-1}TB$ is the pullback bundle, with fibers at $y_1\in N$ given by $(F^{-1}TB)y_1 = T{F(y_1)}B$. The bundle $Hom(TN,F^{-1}TB)$ has a connection induced from the Levi-Civita connection $\nabla^N$ and the pullback connection $\overset{B}{\nabla^F}$. The second fundamental form of $F$ is symmetric and can be expressed as \cite{S7}
	\begin{equation}\label{f}
		(\nabla F_*)(Y_1,Y_2) = \overset{B}{\nabla^F_{Y_1}}F_*Y_2 - F_*(\nabla_{Y_1}^N Y_2)
	\end{equation}
	for all $Y_1,Y_2\in\Gamma(TN)$, where $\overset{B}{\nabla^F_{Y_1}}F_*Y_2\circ F = \nabla^B_{F_*Y_1}F_*Y_2$. In \cite{S7}, \c{S}ahin gave the following Lemma.
	\begin{lemma}
	Let  $F:(N^n,g_N)\rightarrow(B^b,g_B)$ be a map between Riemannian manifolds. Then
	$$g_B((\nabla F_*)(Y_1,Y_2),F_*(Y_3)) = 0,\ \ \ \forall \ Y_1,Y_2,Y_3 \in \Gamma(kerF_*)^\perp.$$
	\end{lemma}
	As a consequence of this Lemma, we get 
	\begin{equation}\label{2.7}
		(\nabla F_{*})(Y_1,Y_2)\in\Gamma (rangeF_*)^\perp.
	\end{equation}
	The tension field \cite{PB} of $F$ is defined as the trace of the second fundamental form of $F$, denoted by $\tau(F)$, and is given by $\tau(F) = \operatorname{trace}(\nabla F_*) = \sum_{i=1}^n (\nabla F_*)(e_i,e_i)$, where $n=\operatorname{dim}(N)$ and ${e_1,e_2,\dots,e_n}$ is the orthonormal frame on $N$ and the bitension field of $F$ \cite{2-har}, denoted by $\tau_2(F)$, is defined as
	\begin{equation}\label{bh}
		\tau_2(F) = -\Delta^F \tau(F) - trace R^B (F_*, \tau(F))F_*.
	\end{equation}
	Moreover, a map $F:(N^n,g_N)\rightarrow(B^b,g_B)$, where $N$ and $B$ are Riemannian manifolds, is biharmonic if and only if the bi-tension field of $F$ vanishes at each point $y_1\in N$.
	\begin{lemma}\label{2.2}\cite{S7}
		Let $F : (N^n, g_N) \rightarrow (B^b, g_B)$ be a Riemannian map between Riemannian manifolds $N$ and $B$. Then, F is umbilical Riemannian map if and only if
		\begin{equation}\label{g}
			(\nabla F_{*})(Y_1,Y_2) = g_N(Y_1,Y_2)H_B
		\end{equation}
		for $Y_1,Y_2 \in \Gamma(rangeF_*)^\perp$ and $H_B$, is nowhere zero, mean curvature vector field on $(rangeF_*)^\perp.$
	\end{lemma}
	
	Let $Y_1$ be a vector field on $N$, and let $Z$ be any section of $(rangeF_*)^\perp$. Then, the orthogonal projection of $\nabla^B_{Y_1}Z$ onto $(rangeF_*)^\perp$, is given by $\nabla^{F\perp}_{Y_1}Z$, where $\nabla^{F\perp}$ is a linear connection on $(rangeF_*)^\perp$ that satisfies $\nabla^{F\perp}g_B=0.$

 Now, for a Riemannian map, we have the following relation \cite{S7}
\begin{equation}\label{h}
	\nabla ^{B}_{F_{*}Y_1}Z=-\mathcal{S}_{Z}F_{*}Y_1+\nabla^{F\perp}_{Y_{1}}Z,
\end{equation}
where $\mathcal{S}_ZF_{*}Y_1$ is the tangential component of $\nabla ^{B}_{F_{*}Y_1}Z$. Thus at $y_1\in N$ we have $\nabla ^{B}_{F_{*}Y_1}Z(y_1)\in T_{F(y_1)}B,~ \mathcal{S}_{Z}F_{*}Y_1(y_1) \in F_{*y_1}(T_{y_1}N)$ and $\nabla^{F\perp}_{Y_{1}}Z(y_1)\in (F_{*y_1}(T_{y_1}N))^\perp.$ Additionally, $\mathcal{S}_{Z}F_{*}Y_1$ is bilinear in $Z$ and $F_{*}Y_1$, and $\mathcal{S}_{Z}F_{*}Y_1$ at $y_1$ depends only on $Z_{y_1}$ and $F_{*y_1}Y_{1y_1}.$ Later, in \cite{M2}, Yadav and Meena stated Clairaut conditions for a geodesic curve on base manifold as follows.

\begin{theorem}\label{thm1}\cite{M2}
	Let $F : (N^n, g_N) \rightarrow (B^b, g_B)$ be a Riemannian map between Riemannian
	manifolds and $\gamma$, $\sigma=F\circ\gamma$ are geodesics on $N$ and $B$, respectively. Then $F$ is Clairaut Riemannian map with $r=e^f$ if and only if any one of the following conditions holds:	
	\begin{enumerate}[label=\roman*)]
		\item $\mathcal{S}_Z F_{*}{Y_1} = -Z(f) F_{*}{Y_1}$, where $F_{*}{Y_1} \in \Gamma(rangeF_*), Z \in \Gamma(rangeF_*)^\perp$ are components
		of $\dot{\sigma}(t).$
		\item $F$ is umbilical map, and has $H_B = -\nabla^B f$, where $f$ is a smooth function on $B$
		and $H_B$ is the mean curvature vector field of $range F_*.$
	\end{enumerate}
\end{theorem} 	

\begin{lemma}\label{2.3}\cite{S7}
		Let $F : (N^n, g_N)\rightarrow (B^b, g_B)$ be a Riemannian map between Riemannian
	manifolds. Then, the tension field $\tau$ of $F$ is
		\begin{equation}\label{tnsn}
			\tau(F)= -n_1F_*(\mu^{\text{ker}F_*}) + n_2H_B,
		\end{equation}
		where $n_1 = dim(kerF_*), \text{and }n_2=rankF,\ \mu^{\text{ker}F_*} \text{ and } H_B$ are the mean curvature vector fields of the distribution $kerF_*$ and  $rangeF_*$, respectively.
	\end{lemma}

	We also recall that the curvature tensor $R^B$ of a Sasakian manifold $(B^b, \varphi, \xi, \eta, g_B)$ with constant $\varphi$-sectional curvature $c_B$ is given by \cite{blairriemannian},
	\begin{equation}\label{21}
		\begin{aligned}
			R^B(B_1, B_2 )B_3 =&\frac{c_B+3}{4}\left\{g_B(B_2, B_3)B_1 - g_B(B_1,B_3)B_2\right\}+\frac{c_B-1}{4}\{\eta(B_1)\eta(B_3)B_2\\&-\eta(B_2)\eta(B_3)B_1 +g_B(B_1,B_3)\eta(B_2)\xi-g_B(B_2,B_3)\eta(B_1)\xi+g_B(\varphi B_2,B_3)\varphi B_1\\&-g_B(B_3,\varphi B_1)\varphi B_2+2g_B(B_1,\varphi B_2)\varphi B_3\},
		\end{aligned}
	\end{equation}
 where $B_1, B_2,B_3\in \Gamma(TB).$

	\section{Clairaut anti-invariant Riemannian maps from Riemannian manifolds to Sasakian manifolds}
	In this section, we define the Clairaut anti-invariant Riemannian map $F$ from a Riemannian manifold $N$ to a Sasakian manifold $B$ and discuss the geometry of such maps. Throughout the section, we consider the Reeb vector field to be in $(range F_*)^\perp$ of $TB$ and $(rangeF_*)^\perp$ as a totally geodesic distribution of $TB$.	
	\begin{definition} 
		\cite{S9} Let $F : (N^n, g_N) \to(B^b, \varphi, \xi, \eta, g_B)$ be a proper Riemannian map from a
		Riemannian manifold $N$ to a Sasakian manifold $B$ with almost contact
		structure $\varphi$. We say that $F$ is an anti-invariant Riemannian map at $y_1 \in N$ if $\varphi(rangeF_{*y_1})\subset(rangeF_{*y_1})^\perp.$ If $F$ is an anti-invariant Riemannian map for every $y_1\in N$ then $F$ is called an anti-invariant Riemannian map.
 \end{definition} 

  In this case, we denote the orthogonal subbundle to $\varphi(rangeF_*)$ in $(rangeF_*)^\perp$ by $\mu$, i.e., $(rangeF_*)^\perp = \varphi(rangeF_*)\oplus \mu$. For any $Z \in\Gamma(rangeF_*)^\perp,$ we have
	\begin{equation}\label{i}
		\varphi Z=\mathcal{B}Z+\mathcal{C}Z,
	\end{equation}
	where $\mathcal{B}Z \in \Gamma(rangeF_*) $ and $\mathcal{C}Z \in \Gamma(\mu$). Moreover, the map $F$ is said to admit horizontal Reeb vector field if $\xi \in (range F_*)^\perp$ and admit vertical Reeb vector field if $\xi \in (range F_*).$
	\begin{definition}\cite{Lag}
		Let $F$ be an anti-invariant map from a Riemannian manifold $(N^n,g_N)$ to an almost contact metric manifold $(B^b, \varphi, \xi, \eta, g_B)$. If $\mu = \{0\}$ or $\mu = \text{span}\{\xi\}$, i.e., $(\text{range}F_*)^{\perp} = \varphi(\text{range}F_*)$ or $(\text{range}F_{*})^{\perp} = \varphi(\text{range}F_*) \oplus\langle \xi \rangle$, respectively, then we call $F$ a Lagrangian Riemannian map.
	\end{definition}

 \begin{definition}
	An antiinvariant Riemannian map $F: (N^n,g_N) \rightarrow (B^b, \varphi, \xi, \eta, g_B)$ from a Riemannian manifold to a Sasakian manifold is called a Clairaut anti-invariant Riemannian map if there is a function $r: B \rightarrow \mathbb{R}^+$ such that for every geodesic $\sigma$ on $B$, the function $(r \circ \sigma) \sin\theta(t)$ is constant, where $F_*Y_1\in \Gamma(rangeF_*)$ and $Z \in \Gamma(rangeF_*)^\perp$ are components of $\dot{\sigma}(t)$, and $\theta(t)$ is the angle between $\dot{\sigma}(t)$ and $Z$ for all $t$.
\end{definition}
	
	\begin{lemma}
		Let $F : (N^n, g_N) \rightarrow (B^b, \varphi, \xi, \eta, g_B)$ be an anti-invariant Riemannian map from a Riemannian manifold $N$ to a Sasakian manifold $B$, and let $\gamma : I \rightarrow N$ be a geodesic on $N$. Then, the curve $\sigma = F \circ \gamma$ is a geodesic on $B$ if and only if
		\begin{equation}\label{j}
         -\mathcal{S}_{\varphi F_{*}Y_1}F_*Y_1-\mathcal{S}_{\mathcal{C}Z}F_*Y_1+\nabla_{Z}^{B}F_*Y + F_*(\mathcal{H}\nabla_{Y_1}^NY)+\eta(Z)F_*Y_1=0,
		\end{equation}
		\begin{equation}\label{k}
			(\nabla F_*)(Y_1,Y)+\nabla^{F\perp}_{Y_1} \varphi F_{*}Y_1 + \nabla^{F\perp}_Z \varphi F_*{Y_1} + \nabla^{F\perp}_Z \mathcal{C}Z + \nabla F^\perp_{Y_1} \mathcal{C}Z-v\xi+\eta(Z)Z = 0,
		\end{equation}
		where $F_*Y_1 \in \Gamma(rangeF_*)$, $Z \in \Gamma(rangeF_*)^\perp$ are components of $\dot{\sigma}(t)$ such that $F_*Y=\mathcal{B}Z$ and $\nabla^B$ is the Levi-Civita connection on $B$, and $\nabla^ {F\perp}$ is a linear connection on $(rangeF_*)^\perp$.
		
		\begin{proof}
			Let $\gamma:I\rightarrow N$ be a geodesic on $N$ and $\sigma=F\circ \gamma$ be a geodesic with speed $\sqrt{v}$ on $B$, where $F_*Y_1 \in \Gamma (rangeF_*)$ and $Z \in \Gamma(rangeF_*)^\perp$ are components of $\dot{\sigma}(t)$. Since $B$ is the Sasakian manifold, we have
			$$\varphi\nabla^B_{\dot{\sigma}}\dot{\sigma}=\nabla^B_{\dot{\sigma}}\varphi\dot{\sigma}-g_B(\dot{\sigma},\dot{\sigma})\xi+\eta(\dot{\sigma})\dot{\sigma}.$$
			Using \eqref{i} and $g_B(\dot{\sigma},\dot{\sigma})=||\dot{\sigma}||^2=v$,
			we get
			$$\varphi\nabla^B_{\dot{\sigma}}\dot{\sigma}=\nabla^B_{F_*Y_1+Z}\varphi (F_*Y_1+Z)-v\xi+\eta(F_*Y_1+Z)(F_*Y_1+Z).$$
			By direct computation, we have
			\begin{equation}\label{l}
					\varphi\nabla^B_{\dot{\sigma}}\dot{\sigma}=\nabla^B_{F_*Y_1}\varphi F_*Y_1+\nabla^B_{F_*Y_1}\varphi Z+\nabla^B_{Z}\varphi F_*Y_1+\nabla^B_{Z}\varphi Z-v\xi+\eta(Z)F_*Y_1+\eta(Z)Z.
			\end{equation}
			Using \eqref{h} and \eqref{i} in \eqref{l}, we get
			\begin{equation}
				\begin{split}\label{m}
					\varphi\nabla^B_{\dot{\sigma}}\dot{\sigma}=&-\mathcal{S}_{\varphi F_*Y_1} F_*Y_1-\mathcal{S}_{\mathcal{C}Z} F_*Y_1+\nabla^B_{F_*Y_1}\mathcal{B}Z+\nabla^{F\perp}_Z\varphi F_*Y_1+\nabla^{F\perp}_{Y_1}\varphi F_*Y_1\\&+\nabla^{F\perp}_{Y_1} \mathcal{C}Z+\nabla^B_{Z} \mathcal{B}Z+\nabla^B_{Z} \mathcal{C}Z-v\xi+\eta(Z)F_*Y_1+\eta(Z)Z.
				\end{split}
			\end{equation}
			We have $g_B(\nabla^B_Z \mathcal{B}Z, V)= 0$, for any $V\in \Gamma(rangeF_*)^\perp$, which implies $\nabla^B_Z \mathcal{B}Z \in \Gamma(rangeF_*)$.
			Let $Y\in \Gamma(kerF_*)^\perp$ such that $F_*Y=\mathcal{B}Z$ and using \eqref{f}, we get 
			\begin{equation}\label{n}
				\nabla^B_{F_*Y_1}\mathcal{B}Z=\overset{B}{\nabla^{F}_{Y_1}}{F_*Y}\circ F=(\nabla F_*)(Y_1,Y)+F_*(\mathcal{H}\nabla^N_{Y_1}Y)
			\end{equation}
			Using \eqref{n} in \eqref{m}, we get
			\begin{equation}\label{o}
				\begin{split}
					\varphi\nabla^B_{\dot{\sigma}}\dot{\sigma}=&-\mathcal{S}_{\varphi F_*Y_1} F_*Y_1-\mathcal{S}_{\mathcal{C}Z} F_*Y_1+(\nabla{F_*})(Y_1,Y)+F_*(\mathcal{H}\nabla^N_{Y_1}Y)+\nabla^{F\perp}_Z\varphi F_*Y_1\\&+\nabla^{F\perp}_{Y_1}\varphi F_*Y_1+\nabla^{F\perp}_{Y_1} \mathcal{C}Z+\nabla^B_{Z} F_*Y+\nabla^B_{Z} \mathcal{C}Z-v\xi+\eta(Z)F_*Y_1+\eta(Z)Z.
				\end{split}
			\end{equation}
			$\text{Now, } \sigma \text{ is geodesic on } B \iff \nabla^B_{\dot{\sigma}}\dot{\sigma}=0 \iff \varphi\nabla^B_{\dot{\sigma}}\dot{\sigma}=0.$
			Therefore, 
			\begin{align*}
				&-\mathcal{S}_{\varphi F_{*}Y_1}F_*Y_1-\mathcal{S}_{\mathcal{C}Z}F_*Y_1 + \nabla_{Z}^{B}F_*Y+F_*(\mathcal{H}\nabla_{Y_1}^NY)+\eta(Z)F_*Y_1+(\nabla F_*)(Y_1,Y)\\& \ \ +\nabla^{F\perp}_{Y_1} \varphi F_{*}Y_1 + \nabla^{F\perp}_Z \varphi F_*{Y_1} + \nabla^{F\perp}_Z \mathcal{C}Z + \nabla F^\perp_{Y_1} \mathcal{C}Z-v\xi+\eta(Z)Z=0,
			\end{align*} which brings the proof to its conclusion.
			
		\end{proof}
	\end{lemma}
	
	\begin{theorem}
		Let $F:(N^n,g_N)\rightarrow (B^b, \varphi, \xi, \eta, g_B)$ be an anti-invariant Riemannian map from a Riemannian manifold $N$ to a Sasakian manifold $B$ and $\gamma, \sigma = F \circ \gamma$ are geodesics on $N$ and $B$, respectively. Then, $F$ is a Clairaut anti-invariant Riemannian map with $r=e^f$ if and only if
		\begin{equation}\label{p}
			\begin{split}
				&g_B(\mathcal{S}_{\varphi F_{*}Y_1}F_*Y_1+\mathcal{S}_{\mathcal{C}Z}F_*Y_1-\eta(Z)F_*Y_1,F_*Y)-g_B((\nabla F_*)(Y_1,Y)+\nabla^{F\perp}_{Y_1} \varphi F_{*}Y_1 \\&+ \nabla^{F\perp}_Z \varphi F_*{Y_1}-v\xi+\eta(Z)Z,\mathcal{C}Z)-g_B(F_*Y_1,F_*Y_1)\frac{d(f\circ\sigma)}{dt}=0,
			\end{split}
		\end{equation} 	
		where $f$ is a smooth function on $B$ and $F_*Y_1\in \Gamma(rangeF_*)$, $Z\in \Gamma(rangeF_*)^\perp$ are the components of $\dot{\sigma}(t)$ and $F_*Y=\mathcal{B}Z$,  for $Y\in \Gamma(kerF_*)^\perp$.
		\begin{proof}
			Let $\gamma: I \rightarrow N$ be a geodesic on $N$ and $\sigma = F \circ \gamma$ be a geodesic with speed $\sqrt{v}$ on $B$ with $F_*Y_1 \in \Gamma(rangeF_*)$ and $Z \in \Gamma(rangeF_*)^\perp$ are components of $\dot{\sigma}(t)$ and $\theta(t)$ denote the angle in $[0, \pi]$ between $\dot{\sigma}$ and $Z$. Therefore, we have	
			\begin{equation}\label{q}
				g_B(Z,Z) = v\cos^2\theta(t),
			\end{equation}
			\begin{equation}\label{r}
				g_B(F_*Y_1,F_*Y_1) = v\sin^2\theta(t).
			\end{equation}
			Differentiating \eqref{q} along $\sigma$, we get
			\begin{equation}\label{s}
				\frac{d}{dt} g_B(Z,Z) = 2g_B(\nabla^B_{\dot{\sigma}} Z, Z) = -2v\sin\theta(t)\cos\theta(t)\frac{d\theta}{dt}.
            \end{equation}	
			Hence, using the Sasakian structure, we get
			$$g_B(\varphi\nabla^B_{\dot{\sigma}} Z, \varphi Z) = -v\sin\theta(t)\cos\theta(t)\frac{d\theta}{dt}.$$
			Now using \eqref{d}, we get
			\begin{equation*}
			g_B(\nabla^B_{\dot{\sigma}}\varphi Z-g_B(\dot{\sigma},Z)\xi+\eta(Z)\dot{\sigma}, \varphi Z) = -v\sin\theta(t)\cos\theta(t)\frac{d\theta}{dt},	
			\end{equation*}
			which implies
			\begin{equation}\label{t}
				g_B(\nabla^B_{\dot{\sigma}}\varphi Z, \varphi Z) = -v\sin\theta(t)\cos\theta(t)\frac{d\theta}{dt}.	
			\end{equation}
			Using \eqref{i} in \eqref{t}, we get
			$$g_B(\nabla^B_{\dot{\sigma}}\mathcal{B} Z,\mathcal{B} Z)+g_B(\nabla^B_{\dot{\sigma}}\mathcal{C} Z,\mathcal{C} Z)= -v\sin\theta(t)\cos\theta(t)\frac{d\theta}{dt}.$$
			Putting $\dot{\sigma}=F_*Y_1+Z$ in the above equation, we get
			\begin{equation*}
   \begin{split}
    &g_B(\nabla^B_{F_*Y_1}\mathcal{B} Z,\mathcal{B}Z)+g_B(\nabla^B_{F_*Y_1}\mathcal{C} Z,\mathcal{C} Z)+g_B(\nabla^B_{Z}\mathcal{B} Z,\mathcal{B}Z)
    \\&+g_B(\nabla^B_{Z}\mathcal{C} Z,\mathcal{C} Z)=-v\sin\theta(t)\cos\theta(t)\frac{d\theta}{dt}. 
   \end{split}  
			\end{equation*}
			Since $(rangeF_*)^\perp$ is totally geodesic, we have
			\begin{equation}\label{u}
				\begin{split}
				 &g_B(\overset{B}{\nabla}{^F_{F_*Y_1}}\mathcal{B} Z\circ F,\mathcal{B}Z)+g_B(\nabla^B_{F_*Y_1}\mathcal{C} Z,\mathcal{C} Z) +g_B(\nabla^B_{Z}\mathcal{B}Z,\mathcal{B}Z)\\&+g_B(\nabla^{F\perp}_{Z}\mathcal{C}Z,\mathcal{C}Z)=-v\sin\theta(t)\cos\theta(t)\frac{d\theta}{dt}.
				\end{split}
			\end{equation}
			Using \eqref{2.7} and \eqref{n} in \eqref{u}, we get
			\begin{equation*}
				g_B(F_*(\mathcal{H}\nabla^N_{Y_1}Y)+\nabla^B_{Z}F_*Y,\: F_*Y)+g_B(\nabla^{F^\perp}_{Y_1}\mathcal{C}Z+\nabla^{F^\perp}_{Z}\mathcal{C}Z ,\: \mathcal{C}Z) =-v\sin\theta(t)\cos\theta(t)\frac{d\theta}{dt}.
			\end{equation*}
			Using \eqref{j} and \eqref{k} in above equation, we obtain	
			\begin{equation}\label{w}
				\begin{split}
					&\ \ g_B(\mathcal{S}_{\varphi F_{*}Y_1}F_*Y_1+\mathcal{S}_{\mathcal{C}Z}F_*Y_1-\eta(Z)F_*Y_1,\:F_*Y)-2g_B((\nabla F_*)(Y_1,Y)\\&+\nabla^{F\perp}_{Y_1} \varphi F_{*}Y_1+ \nabla^{F\perp}_Z \varphi F_*{Y_1}-v\xi+\eta(Z)Z ,\: \mathcal{C}Z) =-v\sin\theta(t)\cos\theta(t)\frac{d\theta}{dt}.
				\end{split}	
			\end{equation}
			Now, $F$ is a Clairaut Riemannian map with $r= e^f$ if and only if $\frac{d}{dt}(e^{f \circ\sigma}\sin \theta) = 0$, \\that is, $$e^{f\circ\sigma}\sin\theta\frac{d(f\circ\sigma)}{dt} + e^{f\circ\sigma}\cos\theta\frac{d\theta}{dt} = 0.$$\\
			Multiplying above equation with non-zero factor $v\sin\theta$ and using \eqref{r}, we get
			\begin{equation}\label{y}
				g_B(F_*Y_1,F_*Y_1)\frac{d(f\circ\sigma)}{dt}=-v\sin\theta(t)\cos\theta(t)\frac{d\theta}{dt}.
			\end{equation}
			Therefore, the assertion \eqref{p} follows from \eqref{w} and \eqref{y}.
		\end{proof}
	\end{theorem}
	
	\begin{theorem}\label{thm2}
		Let $F : (N^n, g_N)\rightarrow (B^b, \varphi, \xi, \eta, g_B)$ be a Clairaut anti-invariant Riemannian map with $r = e^f$ from a Riemannian manifold $N$ to a Sasakian manifold $B$. Then at least one of the following statements is true:
		\begin{enumerate}[label=\roman*)]
			\item $dim(rangeF_*) = 1,$
			\item $f$ is constant on $\varphi(rangeF_*)$, where $f$ is a smooth function on $B$.
		\end{enumerate} 
		\begin{proof}
			Given that $F$ is a Clairaut Riemannian map with $r = e^f$, we can get the following expression by using \hyperref[thm1]{Theorem \ref{thm1}} and \eqref{g} in \eqref{f}:
			\begin{equation}\label{3}
				\overset{B}{\nabla^{F}_{Y_1}}F_{*}Y_2-F_{*}(\nabla_{Y_{1}}^{N}Y_{2}) = -g_N(Y_1,Y_2)\nabla^Bf,
			\end{equation}  
			where $F_*Y_2 \in \Gamma(\text{range}F_*)$ and $Y_1, Y_2 \in \Gamma(\text{ker}F_*)^\perp$. Taking inner product of \eqref{3} with $\varphi F_*X \in \Gamma(\text{range}F_*)^\perp$, we get
			\begin{equation}\label{4}
				g_B(\overset{B}{\nabla^{F}_{Y_1}}F_{*}Y_2 - F_*(\nabla^N_{Y_1} Y_2), \varphi F_*X) = -g_N(Y_1, Y_2)g_B(\nabla^B f, \varphi F_*X).
			\end{equation}     
			Since $\overset{B}{\nabla^F}$ is a Levi-Civita connection, therefore using the metric compatibility condition in \eqref{4}, we get
			\begin{equation}\label{5}
				-g_B(\overset{B}{\nabla^F_{Y_1}}\varphi F_*X, F_*Y_2) = -g_N(Y_1, Y_2)g_B(\nabla^Bf, \varphi F_*X).
			\end{equation} 
			Since $B$ is a Sasakian manifold, therefore using \eqref{d}, we have
			\begin{equation}\label{new}
            \overset{B}{\nabla^F_{Y_1}}\varphi F_*X=g_B(F_*Y_1,F_*X)\xi+\varphi\overset{B}{\nabla^F_{Y_1}}F_*X.
            \end{equation} 
			Using \eqref{new} in \eqref{5}, we get
			$$-g_B(g_B(F_*Y_1,F_*X)\xi+\varphi\overset{B}{\nabla^F_{Y_1}}F_*X, F_*Y_2) = -g_N(Y_1, Y_2)g_B(\nabla^Bf, \varphi F_*X).$$
			By direct computation, we obtain
			\begin{equation}\label{6}
				g_B(\varphi\overset{B}{\nabla^F_{Y_1}}F_*X, F_*Y_2) = g_N(Y_1, Y_2)g_B(\nabla^Bf, \varphi F_*X).
			\end{equation} 
			Using \eqref{c} in \eqref{6}, we get
			$$-g_B(\overset{B}{\nabla^F_{Y_1}}F_*X,\varphi F_*Y_2) = g_N(Y_1, Y_2)g_B(\nabla^Bf, \varphi F_*X).$$
			Using \eqref{3} in the above equation, we find
			\begin{equation}\label{7}
				g_N(Y_1,X)g_B(\nabla^Bf,\varphi F_*Y_2)=g_N(Y_1, Y_2)g_B(\nabla^Bf, \varphi F_*X). 
			\end{equation}
			Now, putting $Y_2=Y_1$ in \eqref{7}, we get
			\begin{equation}\label{8}
				g_N(Y_1,X)g_B(\nabla^Bf,\varphi F_*Y_1)=g_N(Y_1, Y_1)g_B(\nabla^Bf, \varphi F_*X). 
			\end{equation}
			Now interchanging $Y_1$ and $X$ in \eqref{8}, we obtain
			\begin{equation}\label{9}
				g_N(Y_1,X)g_B(\nabla^Bf,\varphi F_*X)=g_N(X,X)g_B(\nabla^Bf, \varphi F_*Y_1). 
			\end{equation}
			From \eqref{8} and \eqref{9}, we get
			$$g_B(\nabla^B f, \varphi F_*Y_1) \left(1 - \frac{g_N(Y_1,Y_1)g_N(X,X)}{g_N(Y_1,X)g_N(Y_1,X)}\right) = 0,$$
			which implies either $dim((kerF_*)^\perp) = 1$ or $g_B(\nabla^B f, \varphi F_*Y_1) = 0$, implies $f$ is constant on $\varphi(rangeF_*)$, which completes the proof.
		\end{proof}
	\end{theorem}
	
	\begin{theorem}
		Let $F : (N^n, g_N)\rightarrow (B^b, \varphi, \xi, \eta, g_B)$ be a Clairaut Lagrangian Riemannian map with $r = e^f$ from a Riemannian manifold $N$ to a Sasakian manifold $B$ such that $dim(rangeF_*) = 1$. Then the following statements are true:
		\begin{enumerate}[label=\roman*)]
			\item $rangeF_*$ is minimal.
			\item $rangeF_*$ is totally geodesic.
		\end{enumerate}
		\begin{proof}
			Since $F$ is a Clairaut Riemannian map with $r = e^f$, from \hyperref[2.2]{Lemma \ref{2.2}}, for any $Y_1 \in (kerF_*)^\perp$, we have
			$$(\nabla F_{*})(Y_1,Y_1) = g_N(Y_1,Y_1)H_B,$$
			where $H_B$ is the mean curvature vector field of $rangeF_*$. Taking the inner product of above equation with $Z \in \Gamma(rangeF_*)^\perp$, we get
			\begin{equation}\label{z}
				g_B(\nabla F_{*})(Y_1,Y_1),Z) = g_N(Y_1,Y_1)g_B(H_B,Z).
			\end{equation}
			Using \eqref{f} in \eqref{z}, we get
			\begin{equation}\label{2}
				g_B(\overset{B}{\nabla^F_{Y_1}} F_{*}Y_1,Z) = g_N(Y_1,Y_1)g_B(H_B,Z).
			\end{equation}
			Using \eqref{b} in \eqref{2}, we get
			\begin{equation}\label{13}
				g(\varphi \overset{B}{\nabla^F_{Y_1}} F_{*}Y_1,\varphi Z)+\eta (\overset{B}{\nabla^F_{Y_1}} F_{*}Y_1)\eta (Z)=g_N(Y_1,Y_1)g_B(H_B,Z),
			\end{equation}	
			Using \eqref{a}, \eqref{c} and \eqref{d} in \eqref{13}, we get
			$$g_B(\overset{B}{\nabla^F_{Y_1}}\varphi F_{*}Y_1,\varphi Z) = g_N(Y_1,Y_1)g_B(H_B,Z).$$
			Since $\overset{B}{\nabla^F}$ is Levi-Civita connection on B. Therefore, using the metric compatibility condition, we get
			\begin{equation}\label{16}
				-g_B(\varphi F_{*}Y_1,\overset{B}{\nabla^F_{Y_1}}\varphi Z) = g_N(Y_1,Y_1)g_B(H_B,Z).
			\end{equation}
			Since, $F$ is a Lagrangian Riremannain map, there exists $Z \in \Gamma(range F_*)^\perp=\Gamma\{\varphi(\text{range}F_*) \oplus\langle \xi \rangle\}$ such that $Z=\varphi F_*X +\xi$ and $\varphi Z=-F_*X$, where $F_*X\in \Gamma (rangeF_*)$, then from \eqref{16}, we get 
			\begin{equation}\label{17}
				-g_B(\varphi F_{*}Y_1,\overset{B}{\nabla^F_{Y_1}}(-F_*X)) = g_N(Y_1,Y_1)g_B(H_B,Z).
			\end{equation}
			Using \eqref{2} in \eqref{17}, we get 
			$$g_B(\varphi F_{*}Y_1,g_N((Y_1,X)H_B)) = g_N(Y_1,Y_1)g_B(H_B,Z).$$	
			Using \hyperref[thm1]{Theorem \ref{thm1}} in the above equation, we get
			$$g_N((Y_1,X))g_B(\varphi F_{*}Y_1,-\nabla^B f) = g_N(Y_1,Y_1)g_B(H_B,Z),$$
			which implies
			\begin{equation}\label{18}
				-g_N(Y_1,X)\varphi F_{*}Y_1(f) = g_N(Y_1,Y_1)g_B(H_B,Z)
			\end{equation}
			Here, $dim(rangeF_*) > 1$, using \hyperref[thm2]{Theorem \ref{thm2}} in \eqref{18}, we get $g_B(H_B, Z) = 0$, which implies $rangeF_*$ is minimal. Also, $H_B=trace(\overset{B}{\nabla^F_{Y_1}}F_*Y_2)$ which shows $rangeF_*$ is totally geodesic.	
		\end{proof}
	\end{theorem}
	
	\begin{theorem}
		Let $F : (N^n, g_N)\rightarrow (B^b, \varphi, \xi, \eta, g_B)$ be a Clairaut anti-invariant Riemannian map from a Riemannian manifold $N$ to a Sasakian manifold $B$ with constant sectional curvature $c_B$. Then $F$ is biharmonic if and only if
		\begin{equation}\label{24}
  \begin{aligned}
			&-n_1{trace}F_*(.)(F_*(\mu^{kerF_*})(f))F_*(.)+n_1{trace}\mathcal{B}\nabla^{F\perp}_{F_*(.)}\nabla^{F\perp}_{F_*(.)} \varphi F_*(\mu^{kerF_*})\\&-n_1{trace}\ \eta(\nabla^{F\perp}_{F_*(.)} \varphi F_*(\mu^{kerF_*}))F_*(.)+n_2{trace}(F_*(.)H_B(f))F_*(.)
                 \\&+n_2{trace}(\nabla^{F\perp}_{F_*(.)} H_B(f))F_*(.)+\left(\frac{c_B+3}{4}\right)n_1(1-n_2)F_*(\mu^{kerF_*})=0
   \end{aligned}
			\end{equation}
   and
		\begin{equation}\label{25}
			\begin{aligned}
				&-n_1(F_*(\mu^{kerF_*})(f))H_B+n_1{trace}({\nabla^{F\perp}_{F_*(.)} \varphi F_*(\mu^{kerF_*})}(f))\varphi F_*(.)\\&+n_1{trace}\mathcal{C}\nabla^{F\perp}_{F_*(.)}\nabla^{F\perp}_{F_*(.)} \varphi F_*(\mu^{kerF_*})+n_2(H_B(f))H_B\\&+n_2{trace}\nabla^{F\perp}_{F_*(.)}\nabla^{F\perp}_{F_*(.)} H_B+\left(\frac{c_B+3}{4}\right)n_2H_B-\left(\frac{c_B-1}{4}\right)\Big\{(n_2)^2\eta(H_B)\xi\\&+n_1{trace}g_N(\mu^{ker{F_*}},.)\xi-3n_2{trace}g_B(H_B,\varphi F_*(.))\varphi F_*(.)\Big\}=0,
			\end{aligned}
		\end{equation}
  where $dim(kerF_*)=n_1$ and $dim((kerF_*)^\perp)=n_2$.
		\begin{proof}
   Let $\{e_i\}_1^{n_1}$ and $\{e_j\}_{1}^{n_2}$ be the orthonormal bases of $(kerF_*)_{y_1}$ and $(kerF_*)_{y_1}^\perp$ at $y_1\in N$. Hence, the Laplacian of $\tau(F)$ is given by
			
			$$\Delta \tau(F)=-\sum_{j=1}^{n_2} \overset{B}{\nabla^F_{e_j}} \overset{B}{\nabla^F_{e_j}} \tau(F).$$
			
			From \hyperref[2.3]{Lemma \ref{2.3}}, we get
			\begin{equation}
				\Delta \tau(F)=-\sum_{j=1}^{n_2} \overset{B}{\nabla^F_{e_j}} \Big \{ \overset{B}{\nabla^F_{e_j}} (-n_1 F_*(\mu^{kerF_*}) + n_2 H_B) \Big\}.	\label{20}
			\end{equation}
			Then, using \eqref{a}, \eqref{d} and \eqref{h} in \eqref{20}, we get
			\begin{equation}\label{30}
   \begin{split}
					\Delta \tau(F) = -\sum_{j=1}^{n_2} &\Big\{n_1\overset{B}{\nabla^F_{e_j}}\varphi (-\mathcal{S}_{\varphi F_*(\mu^{kerF_*})}F_*(e_j)+\nabla^{F\perp}_{F_*(e_j)} \varphi F_*(\mu^{kerF_*}))\\& + n_2\overset{B}{\nabla^F_{e_j}}(-\mathcal{S}_{H_B} F_*(e_j) + \nabla^{F\perp}_{F_*(e_j)} H_B) \Big\},
     \end{split}
			\end{equation}
   where $\{F_*(e_j)\}_{1}^{n_2}$ is an orthonormal basis of of $range F_*$. Since, $F$ is a Clairaut Riemannian map, using \hyperref[thm1]{Theorem \ref{thm1}} and \eqref{a} in \eqref{30}, we get 
   \begin{equation}\label{31}
			  \begin{aligned}
				   \Delta \tau(F) = -\sum_{j=1}^{n_2} &\Big\{n_1\overset{B}{\nabla^F_{e_j}}(-F_*(\mu^{kerF_*})(f)F_*(e_j))
                  +n_1\overset{B}{\nabla^F_{e_j}}\varphi\nabla^{F\perp}_{F_*(e_j)} \varphi F_*(\mu^{kerF_*})\\& 
                  + n_2\overset{B}{\nabla^F_{e_j}}({H_B}(f) F_*(e_j) 
                  + \nabla^{F\perp}_{F_*(e_j)} H_B) \Big\},
            \end{aligned}
	\end{equation}
			Using \eqref{d} and \eqref{h} in \eqref{31}, we get
			\begin{equation}\label{33}
			  \begin{aligned}
				   \Delta \tau(F) =& \ \ \ n_1F_*(\mu^{kerF_*})(f)\sum_{j=1}^{n_2}\overset{B}{\nabla^F_{e_j}}F_*(e_j)+n_1\sum_{j=1}^{n_2}F_*(e_j)(F_*(\mu^{kerF_*})(f))F_*(e_j)
                 \\&-n_1\sum_{j=1}^{n_2}\varphi\Big\{-\mathcal{S}_{\nabla^{F\perp}_{F_*(e_j)} \varphi F_*(\mu^{kerF_*})}F_*(e_j)+\nabla^{F\perp}_{F_*(e_j)}\nabla^{F\perp}_{F_*(e_j)} \varphi F_*(\mu^{kerF_*})\Big\}\\&+n_1\sum_{j=1}^{n_2}\eta(\nabla^{F\perp}_{F_*(e_j)} \varphi F_*(\mu^{kerF_*}))F_*(e_j)-n_2H_B(f)\sum_{j=1}^{n_2}\overset{B}{\nabla^F_{e_j}} F_*(e_j)\\&-n_2\sum_{j=1}^{n_2}(F_*(e_j)H_B(f))F_*(e_j)
                 -n_2\sum_{j=1}^{n_2}\Big\{-\mathcal{S}_{\nabla^{F\perp}_{F_*(e_j)} H_B}F_*(e_j)+\nabla^{F\perp}_{F_*(e_j)}\nabla^{F\perp}_{F_*(e_j)} H_B\Big\}.
            \end{aligned}
	\end{equation}
Again, using \hyperref[thm1]{Theorem \ref{thm1}} and \eqref{i} in \eqref{33}, we get 
 \begin{equation}\label{34}
			  \begin{aligned}
				   \Delta \tau(F) =& \ \ \ n_1(F_*(\mu^{kerF_*})(f))H_B+n_1\sum_{j=1}^{n_2}F_*(e_j)(F_*(\mu^{kerF_*})(f))F_*(e_j)
                 \\&-n_1\sum_{j=1}^{n_2}\Big\{({\nabla^{F\perp}_{F_*(e_j)} \varphi F_*(\mu^{kerF_*})}(f))\varphi F_*(e_j)+\mathcal{B}\nabla^{F\perp}_{F_*(e_j)}\nabla^{F\perp}_{F_*(e_j)} \varphi F_*(\mu^{kerF_*})\Big\}\\&-n_1\sum_{j=1}^{n_2}\Big\{\mathcal{C}\nabla^{F\perp}_{F_*(e_j)}\nabla^{F\perp}_{F_*(e_j)} \varphi F_*(\mu^{kerF_*})-\eta(\nabla^{F\perp}_{F_*(e_j)} \varphi F_*(\mu^{kerF_*}))F_*(e_j)\Big\}\\&-n_2(H_B(f))H_B-n_2\sum_{j=1}^{n_2}(F_*(e_j)H_B(f))F_*(e_j)
                 \\&-n_2\sum_{j=1}^{n_2}\Big\{(\nabla^{F\perp}_{F_*(e_j)} H_B(f))F_*(e_j)+\nabla^{F\perp}_{F_*(e_j)}\nabla^{F\perp}_{F_*(e_j)} H_B\Big\}.
            \end{aligned}
	\end{equation}
			Now, from \eqref{tnsn} and \eqref{21}, we have
			\begin{equation*}
   \begin{aligned}
			\sum_{j=1}^{n_2}R^B(F_*(e_j),\tau(F))F_*(e_j) = &\left(\frac{c_B+3}{4}\right)\Big\{\sum_{j=1}^{n_2}g_B\big(-n_1F_*(\mu^{ker{F_*}}), F_*(e_j)\big)F_*(e_j)\\&-n_2\big(-n_1F_*(\mu^{\text{ker}F_*})+n_2H_B\big)\Big\}+\left(\frac{c_B-1}{4}\right)\Big\{n_2g_B(\xi,n_2H_B)\xi\\&-\sum_{j=1}^{n_2}g_B(-n_1F_*(\mu^{ker{F_*}}),F_*(e_j))\xi\\&+3\sum_{j=1}^{n_2}g_B(n_2\varphi H_B,F_*(e_j))\varphi F_*(e_j)\Big\},
   \end{aligned}
			\end{equation*}
   where $R^B$ is the curvature tensor field of $B$.
   Further simplification of above equation gives
   \begin{equation}\label{12}
   \begin{aligned}
			\sum_{j=1}^{n_2}R^B(F_*(e_j),\tau(F))F_*(e_j) =& -\left(\frac{c_B+3}{4}\right)\Big\{n_1(1-n_2)F_*(\mu^{kerF_*})+n_2H_B\Big\}\\&+\left(\frac{c_B-1}{4}\right)\Big\{(n_2)^2\eta(H_B)\xi+n_1\sum_{j=1}^{n_2}g_N(\mu^{ker{F_*}},e_j)\xi\\&-3n_2\sum_{j=1}^{n_2}g_B\big(H_B,\varphi F_*(e_j)\big)\varphi F_*(e_j)\Big\}.
   \end{aligned}
			\end{equation}
			Therefore, putting \eqref{34} and \eqref{12} in \eqref{bh} and then taking the $rangeF_* $ and $(rangeF_*)^\perp$ parts, we obtain \eqref{24} and \eqref{25}.
		\end{proof}
	\end{theorem}

 \begin{example}
     Let $N=\{(x_1, x_2, x_3)\in \mathbb{R}^3\}$ be a Riemannian manifold with Riemannian metric $$g_N=\begin{pmatrix}
			1 & \frac{1}{\sqrt{2}} & 0 \\
			\frac{1}{\sqrt{2}} & \frac{3}{2}  & 0\\
			0 & 0 & 0  
		\end{pmatrix}$$ and $B=\{(y_1, y_2, y_3)\in \mathbb{R}^3\}$ be a Sasakian manifold with contact structure given as in \hyperref[S]{Example \ref{S}}. Consider a smooth map $F:N\to B$ such that $$F(x_1, x_2, x_3)=\left(0,\frac{2x_1-\sqrt{2}x_2}{\sqrt{3}},\sqrt{3}\right),$$
     after some simple calculations, we get
  $$kerF_*=span\left\{\frac{2e_1+\sqrt{2}e_2}{\sqrt{3}},e_3\right\}$$ and
  $$(kerF_*)^\perp=span\left\{X_1=\frac{2e_1-\sqrt{2}e_2}{\sqrt{3}}\right\},$$
  where $\{e_i\},i=1,2,3$ are standard basis vector fields on $N$. By some computation, we have
  $$rangeF_*=span\left\{F_*X_1={e'_1}=2\frac{\partial}{\partial y_2}\right\},$$
  $$(rangeF_*)^\perp=span\left\{{e'_2}=2\left(\frac{\partial}{\partial y_1}+y_2\frac{\partial}{\partial y_3}\right), e'_3=2\frac{\partial}{\partial y_3}=\xi\right\}$$ and $\varphi e'_1=e'_2$. Here, $g_N(X_1,X_1)=g_B(F_*X_1,F_*X_1)=1.$ Therefore, $F$ is an anti-invariant Riemannian map from the Riemannian manifold to Sasakian manifold.
  In order to show that the defined map is a Clairaut Rriemannian map, we will find a smooth function $f$ satisfying equation $(\nabla F_*)(X_1,X_1)=-g_N(X_1,X_1)\nabla^B f$. Here, it can be seen that $(\nabla F_*)(X_1,X_1)=0$ and $g_N(X_1,X_1)=1$, hence taking $f$ as a constant function, we can say that the given map is a Clairaut anti-invariant Riemannian map.
 \end{example}

 \begin{example}
     Let $N=\{(x_1, x_2, x_3,x_4)\in \mathbb{R}^4\}$ be a Riemannian manifold with Riemannian metric $$g_N=\begin{pmatrix}
			1 & 0 & \frac{-1}{3} & 0 \\
            0 & 1 & 0 & 0 \\
			\frac{-1}{3} & 0 & \frac{2}{3}  & 0\\
			0 & 0 & 0 & 1
		\end{pmatrix}.$$ Let $B=\{(y_1, y_2, y_3,y_4,y_5)\in \mathbb{R}^5\}$ be a Sasakian manifold with contact structure given as in \\ \hyperref[S]{Example \ref{S}}. Consider a smooth map $F:N\to B$ such that $$F(x_1, x_2, x_3,x_4)=\left(0,0,x_1+x_3,x_2+x_4,0\right),$$
     after some simple computations, we obtain
  $$kerF_*=span\{e_1-e_3,e_2-e_4\}$$ and
  $$(kerF_*)^\perp=span\{X_1=e_1+e_3,X_2=e_2+e_4\},$$
  where $\{e_i\},i=1,2,3,4$ are standard basis vector fields on $N$.
  By simple calculations, we get
  $$rangeF_*=span\left\{F_*X_1={e'_1}=2\frac{\partial}{\partial y_3},F_*X_2={e'_2}=2\frac{\partial}{\partial y_4}\right\},$$
  $$(rangeF_*)^\perp=span\left\{{e'_3}=2\left(\frac{\partial}{\partial y_1}+y_3\frac{\partial}{\partial y_5}\right),{e'_4}=2\left(\frac{\partial}{\partial y_2}+y_4\frac{\partial}{\partial y_5}\right),e'_5=2\frac{\partial}{\partial y_5}=\xi\right\}$$ 
  and $\varphi e'_1=e'_3$ and $\varphi e'_2=e'_4$. 
  Here,  $g_N(X_i,X_i)=g_B(F_*X_i,F_*X_i)=1, i=1,2$. Therefore, $F$ is an antiinvariant Riemannian map from the Riemannian manifold to the Sasakian manifold. In order to show that the defined map is a Clairaut Rriemannian map, we will find a smooth function $f$ satisfying equation $(\nabla F_*)(X_1,X_1)=-g_N(X_1,X_1)\nabla^B f$. Here, it can be seen that $(\nabla F_*)(X_1,X_1)=0$ and $g_N(X_1,X_1)=1$, hence taking $f$ as a constant function, we can say that the given map is a Clairaut anti-invariant Riemannian map.
 \end{example}

\section{Acknowledgments}
	The second author is thankful to UGC for providing financial assistance in terms of MANF scholarship vide letter with UGC-Ref. No. 4844/(CSIR-UGC NET JUNE 2019). The third author is thankful to DST Gov. of India for providing financial support in terms of DST-FST label-I grant vide sanction number SR/FST/MS-I/2021/104(C).
 	
	\bibliography{refrns} 

\begin{thebibliography}{}

\bibitem[Allison, 1996]{AD}
Allison, D. (1996).
\newblock Lorentzian {C}lairaut submersions.
\newblock {\em Geometriae Dedicata}, 63:309--319.

\bibitem[Baird and Wood, 2003]{PB}
Baird, P. and Wood, J.~C. (2003).
\newblock {\em Harmonic Morphisms Between {R}iemannian Manifolds}.
\newblock Oxford University Press.

\bibitem[Bishop, 1972]{5}
Bishop, R.~L. (1972).
\newblock {C}lairaut submersions.
\newblock {\em Geometry in Honor of K. Yano, Kinokuniya}, pages 21--31.

\bibitem[Blair, 1976]{B}
Blair, D.~E. (1976).
\newblock {\em Contact Manifolds in {R}iemannian Geometry}.
\newblock Springer Berlin Heidelberg.

\bibitem[Blair et~al., 2010]{blairriemannian}
Blair, D.~E. et~al. (2010).
\newblock {\em Riemannian Geometry of Contact and Symplectic Manifolds}.
\newblock Boston, MA: Birkh{\"a}user Boston: Imprint: Birkh{\"a}user.

\bibitem[Fischer, 1992]{F}
Fischer, A.~E. (1992).
\newblock {R}iemannian maps between {R}iemannian manifolds.
\newblock {\em Contemp. Math}, 132:331--366.

\bibitem[Gray, 1970]{G}
Gray, A. (1970).
\newblock Nearly {K}{\"a}hler manifolds.
\newblock {\em Journal of Differential Geometry}, 4(3):283--309.

\bibitem[Jiang, 1986]{2-har}
Jiang, G.~Y. (1986).
\newblock 2-harmonic maps and their first and second variational formulas.
\newblock {\em Chinese Ann. Math. Ser A}, 7:389--402.

\bibitem[Kumar and Prasad, 2020]{P4}
Kumar, S. and Prasad, R. (2020).
\newblock Semi-slant {R}iemannian maps from cosymplectic manifolds into
  {R}iemannian manifolds.
\newblock {\em Gulf Journal of Mathematics}, 9(1):62--80.

\bibitem[Lee et~al., 2015]{L}
Lee, J., Park, J., {\c{S}}ahin, B., and Song, D.~Y. (2015).
\newblock Einstein conditions for the base space of anti-invariant {R}iemannian
  submersions and {C}lairaut submersions.
\newblock {\em Taiwanese Journal of Mathematics}, 19(4):1145--1160.

\bibitem[Li et~al., 2022]{YL}
Li, Y., Prasad, R., Haseeb, A., Kumar, S., and Kumar, S. (2022).
\newblock A study of {C}lairaut semi-invariant {R}iemannian maps from
  cosymplectic manifolds.
\newblock {\em Axioms}, 11(10):503.

\bibitem[Meena and Yadav, 2023]{M2}
Meena, K. and Yadav, A. (2023).
\newblock {C}lairaut {R}iemannian maps.
\newblock {\em Turkish Journal of Mathematics}, 47(2):794--815.

\bibitem[O'Neill, 1966]{O}
O'Neill, B. (1966).
\newblock The fundamental equations of a submersion.
\newblock {\em Michigan Mathematical Journal}, 13(4):459--469.

\bibitem[Pastore et~al., 2004]{8}
Pastore, A.~M., Falcitelli, M., and Ianus, S. (2004).
\newblock {\em {R}iemannian Submersions and Related Topics}.
\newblock World Scientific.

\bibitem[Prasad and Kumar, 2017]{P2}
Prasad, R. and Kumar, S. (2017).
\newblock Slant {R}iemannian maps from {K}enmotsu manifolds into {R}iemannian
  manifolds.
\newblock {\em Global J. Pure Appl. Math}, 13:1143--1155.

\bibitem[Prasad and Kumar, 2018]{P3}
Prasad, R. and Kumar, S. (2018).
\newblock Semi-slant {R}iemannian maps from almost contact metric manifolds
  into {R}iemannian manifolds.
\newblock {\em Tbilisi Mathematical Journal}, 11(4):19--34.

\bibitem[Prasad and Pandey, 2017]{P1}
Prasad, R. and Pandey, S. (2017).
\newblock Slant {R}iemannian maps from an almost contact manifold.
\newblock {\em Filomat}, 31(13):3999--4007.

\bibitem[{\c{S}}ahin, 2010]{S9}
{\c{S}}ahin, B. (2010).
\newblock Invariant and anti-invariant {R}iemannian maps to {K}{\"a}hler
  manifolds.
\newblock {\em International Journal of Geometric Methods in Modern Physics},
  7(03):337--355.

\bibitem[{\c{S}}ahin, 2017a]{S11}
{\c{S}}ahin, B. (2017a).
\newblock Circles along a {R}iemannian map and {C}lairaut {R}iemannian maps.
\newblock {\em Bulletin of the Korean Mathematical Society}, 54(1):253--264.

\bibitem[{\c{S}}ahin, 2017b]{S7}
{\c{S}}ahin, B. (2017b).
\newblock {\em {R}iemannian submersions, {R}iemannian maps in {H}ermitian
  geometry, and their applications}.
\newblock Academic Press.

\bibitem[Ta{\c{s}}tan, 2017]{Lag}
Ta{\c{s}}tan, H.~M. (2017).
\newblock Lagrangian submersions from normal almost contact manifolds.
\newblock {\em Filomat}, 31(12):3885--3895.

\bibitem[Taṣtan and Gerdan, 2017]{HMT}
Taṣtan, H.~M. and Gerdan, S. (2017).
\newblock {C}lairaut anti-invariant submersions from {S}asakian and {K}enmotsu
  manifolds.
\newblock {\em Mediterranean Journal of Mathematics}, 14:1--17.

\bibitem[Yadav and Meena, 2022a]{AY2}
Yadav, A. and Meena, K. (2022a).
\newblock {C}lairaut anti-invariant {R}iemannian maps from {K}{\"a}hler
  manifolds.
\newblock {\em Mediterranean Journal of Mathematics}, 19(3):97.

\bibitem[Yadav and Meena, 2022b]{AY}
Yadav, A. and Meena, K. (2022b).
\newblock {C}lairaut invariant {R}iemannian maps with {K}{\"a}hler structure.
\newblock {\em Turkish Journal of Mathematics}, 46(3):1020--1035.

\bibitem[Zaidi et~al., 2023]{Azaidi}
Zaidi, A., Shanker, G., and Yadav, A. (2023).
\newblock Conformal anti-invariant riemannian maps from or to sasakian
  manifolds.
\newblock {\em Lobachevskii Journal of Mathematics}, 44(4):1518--1527.

\end{thebibliography}
	\bibliographystyle{apalike}
	
\end{document}